\documentclass[11pt, a4paper]{amsart}

\usepackage{amsmath}
\usepackage{amsthm}
\usepackage{amssymb}
\usepackage{srcltx}
\usepackage{url}
\usepackage[hidelinks]{hyperref}
\usepackage{tikz}

\newtheorem{lemma}{Lemma}[section]
\newtheorem{thm}[lemma]{Theorem}

\theoremstyle{definition}

\newtheorem{rem}[lemma]{Remark}

\numberwithin{equation}{section}

\begin{document}

\title{Model Predictive Control for Regular Linear Systems} 

\author{Stevan Dubljevic}
\address{University of Alberta, Chemicals and Materials
  Engineering Department , Edmonton, AB T6G 2V4, Canada}
\email{stevan.dubljevic@ualberta.ca}
\author{Jukka-Pekka Humaloja}
\address{ampere University, Computing Sciences, Mathematics, P.O. Box
  692, 33014 Tampere University, Finland}
\email{jukka-pekka.humaloja@tuni.fi}

\thispagestyle{plain}

\begin{abstract}                         
The present work extends known finite-dimensional constrained optimal control realizations to the realm of well-posed regular linear infinite-dimensional systems modelled by partial differential equations. The structure-preserving Cayley-Tustin transformation is utilized to approximate the continuous-time system by a discrete-time model representation without using any spatial discretization or model reduction. The discrete-time model is utilized in the design of model predictive controller accounting for optimality, stabilization, and input and output/state constraints in an explicit way. The proposed model predictive controller is dual-mode in the sense that predictive controller steers the state to a set where exponentially stabilizing unconstrained feedback can be utilized without violating the constraints. The construction of the model predictive controller leads to a finite-dimensional constrained quadratic optimization problem easily solvable by standard numerical methods. Two representative examples of partial differential equations are considered.
\end{abstract}

\subjclass[2010]{49K20, 93C05 (93B17, 49N10, 35L05)}

\keywords{infinite-dimensional systems, modeling and control optimization, controller constraints and structure, model predictive control, regular linear systems, Cayley-Tustin transform}

\thanks{This work was initiated while J.-P. Humaloja
  was visiting University of Alberta in 2017 and in 2018. The first visit was funded by
  the Doctoral Program of Engineering and Natural
  Sciences of Tampere University of Technology (TUT) and the second one was
  supported by the International HR services of TUT. The research is
  supported by the Academy of Finland Grant  number 310489 held by Lassi Paunonen.}

\maketitle

\section{Introduction} \label{sec:intro}

The concept of \textit{regular linear systems} came about at the turn
of 1990's by the work of George Weiss \cite{Wei89b, Wei94a,
  Wei94b}. This subclass of abstract linear systems is essentially the
Hilbert space counterpart of the finite-dimensional systems described
by the state-space equations:
\begin{subequations}
\begin{align}
\dot{x}(t) & = Ax(t) + Bu(t), \quad x(0) = x_0 \\
y(t) & = Cx(t) + Du(t),
\end{align}
\label{eq:linsys}%
\end{subequations}
where, however, the operators $A$, $B$ and $C$ may be
unbounded. Regular linear systems are often encountered in the study
of partial differential equations (PDEs) with boundary controls and boundary
observations, and they cover a large class of abstract systems of
practical interest.

Over the past decade, there have been several attempts to address control of distributed parameter systems within an input and/or state constrained optimal control setting.  There are several works on dynamical analysis and optimal control of hyperbolic PDEs, most notably the work of Aksikas et al. on optimal linear quadratic feedback controller design for hyperbolic DPS \cite{aksikas2009lq,aksikas2008optimal,Guay.2004}. Other contributions considered optimal and model predictive control applied to Riesz spectral systems (parabolic and higher order dissipative PDEs) with a separable eigenspectrum of the underlying dissipative spectral operator and successfully designed algorithms that account for the input and state constraints \cite{Yang2013136225,Liu20141671,Dubljevic2006IJRNC}. In prior contributions, some type of spatial approximation is applied to the PDE models to arrive at finite-dimensional models utilized in the controller design. As it will be claimed and demonstrated in the subsequent sections, the linear distributed parameter system can be treated intact and controller design can be accomplished without any spatial model approximation or reduction.

The research area of model predictive control (MPC) and contributions associated with this design methodology has flourished over past two decades \cite{rawlings2017model,zbMATH05866352,MAYNE2000789,zbMATH06431684}. The appealing nature of applying to the state the first control input in a finite sequence of control inputs obtained as a solution of an online constrained, discrete-time, optimal control problem with explicit account for the control and state constraints, and achieving stability by adding a terminal cost or terminal constraints, or by extending the horizon of the optimal control problem, is well understood and explored \cite{MAYNE2000789,zbMATH06431684,zbMATH00508571} but could not be easily extended to the DPS setting. Contributions on this topic include, e.g., \cite{ItoKun02} where the terminal penalty approach was analyzed for infinite-dimensional systems with bounded controls, \cite{PhaGeo14} where MPC was formulated for boundary controlled hyperbolic systems with zero terminal constraint, and \cite{Phageo12} where MPC was studied for scalar nonlinear hyperbolic systems. For contributions on MPC for parabolic PDEs, see e.g. \cite{LaoEll14} and the references therein.

In addition to the aforementioned contributions, unconstrained nonlinear MPC for finite- and infinite-dimensional systems has been considered in \cite{zbMATH05702475} with emphasis on the computational complexity of the optimization problem. However, the clear link between the discrete constrained optimization based MPC design, the well-understood modelling of distributed parameter systems described by PDEs, and the well-established control theory of linear DPS has not yet been established apart from the recent work by the authors \cite{HumDubACC18, XuQDub17}.

Motivated by the preceding, in this contribution, the model predictive control for regular linear systems is developed. In particular, the essential feature of the discrete-time infinite-dimensional representation necessary in the MPC design preserving the continuous-time system properties is established by applying the Caley-Tustin (CT) \cite{HavMal07} time discretization, implying that no spatial discretization or model reduction is required. At the core of the  CT transformation, one can find the application of a Crank-Nicolson type time discretization scheme which is a well-know implicit midpoint integration rule that is symmetric, symplectic (Hamiltonian preserving) \cite{hairer2002geometric}, and guarantees structure preserving numerical integration so that stability and controllability are not altered by the discrete-time infinite-dimensional model representation \cite{CurOos97}. Furthermore, boundary and/or point actuation transformed to the discrete-setting yields bounded operators. Due to the input/output convergence and the stability preserving properties of the CT transformation \cite{HavMal07}, the developed discrete-time control laws can be applied to achieve stabilization also for the original continuous-time systems. However, optimality of the controls in continuous time cannot be guaranteed.

As the first main contribution, we develop a linear model predictive control algorithm for stable regular linear systems and prove optimality and stability of the proposed design (Theorem \ref{thm:1}). The design is demonstrated on a numerical example of the one-dimensional wave equation. The proposed design was introduced in \cite{XuQDub17} where it was applied to specific PDE models. Here we extend the design to the broad class of regular linear systems and provide sufficient general conditions under which we prove stability and optimality of the proposed design, in which way we significantly extend the results of \cite{XuQDub17}.

As the second main contribution, an MPC-based control design is presented to achieve constrained stabilization of exponentially stabilizable regular linear systems with proof on optimality and stability (Theorem \ref{thm:2}). The proposed design belongs to the class of dual-mode control \cite{262032,751369} implying that the model predictive controller steers the state to the neighborhood of the origin where local unconstrained stabilizing feedback can be applied without violating the input constraints. A stabilizing terminal penalty is added to the MPC formulation to guarantee stabilizability while no terminal constraints are imposed. Stabilization of parabolic systems was considered in the MPC setting in \cite{XuQDub17}, but there the unstable nodes were tackled with an additional equality constraint, which is only applicable for a finite number of unstable eigenvalues. Our approach here is different as it utilizes an additional terminal penalty function and can be applied to arbitrary exponentially stabilizable regular linear systems, which is a major improvement to  \cite{XuQDub17} in terms of the broad class of systems the design can be applied to. The proposed control design is demonstrated on a simulation study of a tubular reactor that has infinitely many unstable eigenvalues.

The structure of the paper is as follows. In Section \ref{sec:math}, we present the notation, the mathematical preliminaries concerning regular linear systems and the Cayley-Tustin time discretization
scheme. In Section \ref{sec:mpc}, we present the MPC problem, and in Sections \ref{mpc:stab} and \ref{mpc:ustab}, stability and optimality results of the proposed MPC and dual-mode control designs are presented. In Section \ref{sec:waveeq}, we present, as an example of a stable system, the wave equation on a one-dimensional spatial
domain and compute the operators corresponding to the Cayley-Tustin transform. Furthermore, in Section \ref{waveeq:lyap}, we derive a
solution of the Lyapunov equation for the wave equation as required by the proposed MPC design. The performance of the MPC is demonstrated by
numerical simulations of the controlled wave equation in Section \ref{waveeq:sim}. In Section \ref{sec:tub}, the dual-mode controller design is
demonstrated on an unstable tubular reactor which is successfully stabilized by the proposed control strategy. Finally, the paper is concluded in Section
\ref{sec:concl}. 

\section{Mathematical Preliminaries} \label{sec:math}

\subsection{Notation} \label{sec:not}

Here $\mathcal{L}(X,Y)$ denotes the set of bounded linear operators
from the normed space $X$ to the normed space $Y$. The domain, range,
kernel and resolvent of a linear operator $A$ are denoted by
$\mathcal{D}(A)$, $\mathcal{R}(A)$, $\mathcal{N}(A)$ and $\rho(A)$,
respectively. For a linear operator $A: \mathcal{D}(A) \subset X \to
X$ and a fixed $s_0 \in \rho(A)$, define the scale spaces $X_1 :=
(\mathcal{D}(A), \| (s_0 - A)\cdot\|)$ and $X_{-1} = \overline{(X,
  \|(s_0-A_{-1})^{-1}\cdot\|)}$, where $A_{-1}$ is the extension of
$A$ to $X_{-1}$ \cite[Sec. 2.10]{TucWeiBook}. The scale
spaces are related by $X_1 \subset X \subset X_{-1}$ where the
inclusions are dense and with continuous embeddings. The $\Lambda$-extension of an operator $P$ is denoted by $P_{\Lambda}$ (see \eqref{eq:clam}).

\subsection{Regular Linear Systems} \label{sec:rls}

Consider a well-posed linear system $(A,B,$ $C,D)$, where 
$A: \mathcal{D}(A) \subset X \to X$ is the generator of a
$C_0$-semigroup, $B\in \mathcal{L}(U,X_{-1})$ is the control operator, $C \in
\mathcal{L}(X_1,Y)$ is the observation operator, and $D \in
\mathcal{L}(U,Y)$. We assume that the spaces $X$, $U$,
and $Y$ are separable Hilbert spaces and that $U$ and $Y$ are finite-dimensional.

The operator $B$ is called an \textit{admissible input operator} for $A$
if for some $\tau > 0$, the operator $\Phi_{\tau}\in
\mathcal{L}(L^2(0,\infty;U), X_{-1})$ defined as \cite[Sec. 4.2]{TucWeiBook}:
$$
\Phi_{\tau}u = \int\limits_0^{\tau}T(\tau - s)Bu(s)ds,
$$
satisfies $\mathcal{R}(\Phi_{\tau}) \subset X$. Correspondingly, the
operator $C$ is called an \textit{admissible output operator} for $A$
if for some $\tau > 0$, there exists a $K_{\tau}$ such that \cite[Sec. 4.3]{TucWeiBook}:
$$
\int\limits_0^{\tau}\| CT(s)x \|^2ds \leq K_{\tau}\| x
\|^2, \qquad \forall x \in \mathcal{D}(A).
$$
Furthermore, if there exists a $K$ such that $K_{\tau} \leq K$ for
all $\tau > 0$, then $C$ is called \textit{infinite-time
  admissible}. Infinite-time admissibility of the observation operator is required later on for the solvability of the Lyapunov equation \eqref{eq:lyacont}.
  
Let $G$ denote the transfer function of the system $(A,B,C,D)$. The
transfer function is called \textit{regular} if $\displaystyle \lim_{\lambda \to \infty}G(\lambda)u = Du$ for all $u \in U$ \cite[Thm. 1.3]{Wei94b}, in which case $(A,B,C,D)$ is called a \textit{regular linear
  system}. The transfer function $G$ of a regular system is given by:
$$
G(s) := C_{\Lambda}(s - A_{-1})^{-1}B + D,
$$
where $C_\Lambda$ denotes the \textit{$\Lambda$-extension} of the operator $C$
defined as \cite{Wei94a}:
\begin{equation}
\label{eq:clam}
C_{\Lambda}x = \lim_{\lambda \to \infty}\lambda C(\lambda - A)^{-1}x,
\end{equation}
the domain of which consists of those elements $x \in X$
for which the limit exists. Regular linear system have the usual state space presentation \eqref{eq:linsys}, where $C$ must be replaced by $C_\Lambda$. Throughout this paper, we assume that we are dealing with regular
linear systems with admissible $B$ and $C$. However, we note that the  admissibility of $B$ is only needed for the state-feedback stabilization in Section \ref{mpc:ustab} and can be lifted elsewhere.

\subsection{Cayley-Tustin Time Discretization} \label{sec:ct}

Consider a system given in \eqref{eq:linsys}. Given a time
discretization parameter $h>0$, the Tustin time discretization of
\eqref{eq:linsys} is given by
\begin{align*}
\frac{x(jh) - x((j-1)h)}{h} & \approx A \frac{x(jh) -
                                          x((j-1)h)}{2} + Bu(jh) \\
  y(jh) & \approx C \frac{x(jh) - x((j-1)h)}{2} + Du(jh)
\end{align*}
for $j\geq 1$, where we omitted the spatial dependence of $x$ for brevity.
Let $u_j^{(h)}/\sqrt{h}$ be the approximation of
$u(t)$ on the interval $t\in ((j-1)h,jh)$, e.g., by the
mean value sampling used in \cite{HavMal07}: 
$$
\frac{u_j^{(h)}}{\sqrt{h}} = \frac{1}{h}\int\limits_{(j-1)h}^{jh}u(t)dt.
$$
It has been shown in \cite{HavMal07} that the Cayley-Tustin
discretization is a convergent time discretization scheme for
input-output stable system nodes satisfying $\dim U = \dim Y = 1$ in
the sense that $y_j^{(h)}/\sqrt{h}$ converges to $y(t)$ in several
different ways as $h\to 0$. The discussion in \cite[Sec. 6]{HavMal07} further
implies that the same holds for any finite dimensional $U$ and
$Y$. Thus, writing $y_j^{(h)}/\sqrt{h}$ and
$u_j^{(h)}/\sqrt{h}$ in place of $y(jh)$ and $u(jh)$, respectively,
simple computations yield the \textit{Cayley-Tustin
  discretization} of \eqref{eq:linsys} as:
\begin{align*}
  x(k) & = A_dx(k-1) + B_du(k), \quad x(\zeta,0)=x_0(\zeta) \\
  y(k) & = C_dx(k-1) + D_du(k),
\end{align*}
where:
$$
\begin{bmatrix}
A_d & B_d \\ C_d & D_d
\end{bmatrix}
   := \begin{bmatrix}
   -I + 2\delta(\delta - A)^{-1} & \sqrt{2\delta}(\delta - A_{-1})^{-1}B \\
  \sqrt{2\delta}C(\delta - A)^{-1} & G(\delta)
\end{bmatrix}
$$
and $\delta := 2/h$. Clearly one must have $\delta \in \rho(A)$, so that
the resolvent operator is well-defined. Thus, for a large enough $\delta$, the discretization can be applied to unstable systems as well.

When a discrete-time control law $u(k)$ is obtained, it can be transferred to continuous time, e.g., by defining $u(t) = u(k)/\sqrt{h}$ for $t \in [kh, (k+1)h]$, and based on the input/output convergence of the CT discretization \cite{HavMal07}, we have that under the control $u(t)$, the continuous-time output behaves approximately as $y(t) \approx y(k)/\sqrt{h}$. That is, the continuous-time output approximately follows the discrete-time output under the discrete-time control law. As the CT transform also preserves input/output stability \cite{HavMal07}, a stabilizing discrete-time control law $u(k)$ designed for the CT discretized system can be used to stabilize the original continuous-time as well.

\section{Model Predictive Control} \label{sec:mpc}

The moving horizon regulator is based on a similar formulation
emerging from the finite-dimensional system theory (see e.g.
\cite{MusRaw93}). A corresponding controller in the infinite-dimensional
case is presented, e.g., in \cite{XuQDub17}. At a given sampling time
$k$, the objective function with constraints is given by:
\begin{equation}
  \label{eq:mpc}
  \begin{aligned}
 \min_{u} \sum_{j=k+1}^{\infty}& \langle y_{k+j}, Qy_{k+j} \rangle_{Y} + \langle
                 u_{k+j}, Ru_{k+j} \rangle_{U} \\
  \mbox{s.t.} \quad x_{j}& = A_dx_{j-1} + B_du_{j} \\
  y_{j} & = C_dx_{j-1} + D_du_{j} \\
  u_{\min}& \leq u_{j} \leq u_{\max}, \qquad 
  y_{\min} \leq y_{j} \leq y_{\max},
\end{aligned}
\end{equation}
where $Q$ and $R$ are positive self-adjoint weights on the
outputs $y_{j}$ and inputs $u_{j}$, respectively. Here it is assumed
for simplicity that $U$ and $Y$ are (finite-dimensional) real-valued spaces. For
consideration of the MPC with complex input and output spaces, see
\cite{HumDubACC18}, where the authors considered MPC for the
Schr\"odinger equation.

The infinite-horizon objective function \eqref{eq:mpc} can be cast
into a finite-horizon  objective function under certain assumptions on
the inputs beyond the control horizon. Furthermore, a penalty term
needs to be added to the objective function to account for the inputs
and outputs beyond the horizon. We will present two approaches on this
depending on the stability of the original plant.

\subsection{Stable systems} \label{mpc:stab}

If $A$ is the generator of a (strongly) stable
$C_0$-semi- group, we may assume that the input is zero beyond the
control horizon $N$, i.e., $u_{k+N+i} = 0, \forall i\in \mathbb{N}$, and
add a corresponding output penalty term. Under the
assumption that $C$ is infinite-time admissible for $A$,
the terminal output penalty term can be
written as a state penalty term, so that the finite-horizon objective
function is given by:
\begin{equation}
  \label{eq:mpcfin}  
  \min_{u^N} \sum_{j=k+1}^{k+N} \langle y_{j}, Qy_{j}
  \rangle_{Y}  + \langle
                 u_{j}, Ru_{j} \rangle_{U}  + \langle x_{k+N},
                 \bar{Q}x_{k+N} \rangle_X
\end{equation}
with the same constraints as in \eqref{eq:mpc}, and where $N$ is the length of the control horizon.

The operator $\bar{Q}$ can be calculated from the positive self-adjoint
solution of the following discrete-time Lyapunov equation:
\begin{equation}
  \label{eq:lyadisc}
  A_d^{*} \bar{Q}A_d - \bar{Q} = -C_d^{*} QC_{d},
\end{equation}
or equivalently (see \cite[Thm. 2.4]{CurOos97}) the continuous-time Lyapunov equation:
\begin{equation}
  \label{eq:lyacont}
  A^{*} \bar{Q} + \bar{Q}A = -C^{*}QC
\end{equation}
on the dual space of $X_{-1}$. By the assumed infinite-time admissibility of $C$
and the stability of $A$, the operator $\bar{Q} \in \mathcal{L}(X)$ given by \cite[Thm. 5.1.1]{TucWeiBook}:
\begin{equation}
  \label{eq:lyasol}
  \bar{Q}x = \lim_{\tau \to \infty}
  \int\limits_0^{\tau}T^{*}(t)C^{*}QCT(t)xdt, \quad \forall x \in \mathcal{D}(A),
\end{equation}
is the unique positive self-adjoint solution of the Lyapunov equations
\eqref{eq:lyacont} and \eqref{eq:lyadisc}.

Now that we have established that the finite-horizon objective
function \eqref{eq:mpcfin} is well-defined, to further manipulate the
objective function \eqref{eq:mpcfin} we introduce the notation $Y_k := (y_{k+n})_{n=1}^N \in Y^N$ and $U_k := (u_{k+n})_{n=1}^N \in U^N$. Hence, a manipulation of the objective function \eqref{eq:mpcfin} leads to the following quadratic optimization problem:
\begin{equation}
  \label{eq:mpcmat}
  \min_{U_{k}} \, \langle U_k, HU_k \rangle_{U^{N}} + 2\langle U_k, Px_k
  \rangle_{U^{N}} + \langle x_k, \bar{Q}x_k \rangle_X,
\end{equation}
where $H \in \mathcal{L}(U^{N})$ is positive and self-adjoint given by:
$$
h_{i,j} = \begin{cases}
		D_d^*QD_d + B_d^*\bar{Q}B_d + R, & \mbox{for } i = j \\
		D_d^*QC_dA_d^{i-j-1}B_d + B_d^*\bar{Q}A_d^{i-j}B_d, & \mbox{for } i > j \\
		h_{j,i}^*, & \mbox{for } i < j
	\end{cases}
$$
and $P \in \mathcal{L}(X, U^{N})$ is given by $P = (D_d^*QC_dA_d^{k-1} + B_d^*\bar{Q}A_d^k)_{k=1}^N$

The objective function \eqref{eq:mpcmat} is subjected to constraints $U_{\min} \leq U_k \leq U_{\max}$ and $Y_{\min} \leq \left(SU_k + Tx_k\right) \leq Y_{\max}$
which can be written in the form:
\begin{equation}
  \label{eq:linconstr}
\begin{bmatrix}
I\\  -I \\  S \\ - S
\end{bmatrix} U_k \leq \begin{bmatrix}
U_{\max}  \\  -U_{\min}\\ Y_{\max} - Tx_k \\ -Y_{\min} + Tx_k
\end{bmatrix},
\end{equation}
where $S \in \mathcal{L}(U^{N}, Y^N)$ is given by:
$$
	s_{i,j} = \begin{cases}
		D_d, & \mbox{for } i = j \\
		C_dA_d^{i-j-1}B_d, & \mbox{for } i > j \\
		0, & \mbox{for } i < j
	\end{cases}
$$
and $T \in \mathcal{L}(X, Y^{N})$ is given by $T = (C_dA_d^{k-1})_{k=1}^N$.

Considering a finite-dimensional input space $U = \mathbb{R}^m$, the
inner products in the objective function given in
\eqref{eq:mpcmat} are simply vector products, and we have a finite
dimensional quadratic optimization problem:
\begin{equation}
\label{eq:jkfin}
\min_{U_k} J(U_k,x_k) =  U_k^THU_k + 2U_k^T(Px_k).
\end{equation}
Note that the term $\left\langle x_k, \bar{Q}x_k \right\rangle_X$ can
be neglected as $x_k$ is the initial condition for step $k+1$ and
cannot be affected by the control input. Furthermore, as all the
operators related to the objective function and the linear constraints
are bounded under the standing assumptions, the quadratic optimization
problem is exactly of the same
form as the ones obtained for finite-dimensional systems. Thus, we
obtain the convergence and stability results for free by the MPC
theory on finite-dimensional systems (see e.g.
\cite{751369}). To highlight this observation, we
present the following result:

\begin{thm}
  \label{thm:1}
Assume that $A$ is the
generator of a strongly stable $C_0$-semigroup and that $C$ is an
infinite-time admissible observation operator for $A$. Then, the
input sequence $(U_k)$ (and hence the sequence $(u_k)$)
obtained as the solution of the feasible quadratic optimization problem \eqref{eq:jkfin}
with constraints \eqref{eq:linconstr} converges to zero.
\begin{proof}
  By the preceding argumentation, the resulting MPC problem is equivalent
  to a finite-dimensional one, and thus, the result follows from standard finite-dimensional MPC theory.
\end{proof}
\end{thm}

\begin{rem}
Due to the assumed strong stability of the semigroup generated by $A$,
the state of the system under the MPC control law goes asymptotically
to zero for all initial states $x_0 \in \mathcal{D}(A)$ for which the problem is feasible as the
 control inputs decay to zero by Theorem \ref{thm:1}.
\end{rem}

\subsection{Exponentially stabilizable systems} \label{mpc:ustab}

Let us now assume that the pair $(A,B)$ is exponentially stablizable,
i.e., there exists an admissible feedback operator $K \in
\mathcal{L}(X_1,U)$ such that $A + BK_{\Lambda}$ is the generator of
an exponentially stable $C_0$-semigroup
\cite[Def. 3.1]{WeiCur97}, and moreover that $C$ is admissible for $A+BK_\Lambda$. As here input and output spaces are finite-dimensional, the considered systems are uniformly line-regular in the sense of \cite[Def. 6.2.3]{MikPhd}, and thus, the optimal (in terms of minimizing the continuous version of \eqref{eq:mpc}) stabilizing state feedback operator is
obtained using the maximal solution $\bar{R} \in \mathcal{L}(X)$ of the continuous-time Riccati equation \cite[Def. 10.1.2]{MikPhd}:
\begin{equation}
  \label{eq:cric}
  K^{*}SK = A^{*} \bar{R} + \bar{R}A + C^{*}QC
\end{equation}
on $\mathcal{D}(A)$, where $S := R + D^{*}QD$ and $K := -S^{-1}(B_{\Lambda}^{*} \bar{R} +
D^{*}QC)$ yields the optimal feedback operator. Moreover, by \cite[Sect. 3]{CurOos97} the solutions of \eqref{eq:cric} are
equivalent to the solutions of the discrete-time Riccati equation:
\begin{equation}
  \label{eq:dric}
   K_d^{*}S_dK_d = A_d^{*} \bar{R}A_d - \bar{R} + C_d^{*}QC_d,
\end{equation}
where $S_d := B_d^{*} \bar{R}B_d + R + D_d^{*}QD_{d}$ and $K_d :=
-S_d^{-1}(A_d \bar{R}B_d + D_d^{*}QC_d)$ yields the optimal state feedback
for the discrete-time system with the maximal $\bar{R}$ so that $A_d+B_dK_d$ corresponds to
the Cayley-Tustin discretization of $A+BK_{\Lambda}$. Thus, as $K$ is
a stabilizing feedback for $(A,B)$, equivalently
$A_{K_d} := A_d+B_dK_d$ is asymptotically stable \cite[Lem 2.2]{CurOos97}.

Returning to the MPC problem, we assume that the optimal state feedback is
utilized beyond the control horizon, i.e., $u_{k+N+i} =
K_dx_{k+N+i-1}, \forall i\in \mathbb{N}$. Thus, 
the input and output terminal penalties can be expressed as state
terminal penalties by solving the discrete-time Lyapunov equations:
\begin{align*}
  A_{K_d}^{*} \bar{Q}_1 A_{K_d} - \bar{Q}_1 & = -K_d^{*}RK_d \\
  A_{K_d}^{*} \bar{Q}_2A_{K_d} - \bar{Q}_2 & = -(C_d + K_dD_d)^{*}Q (C_d+K_dD_d)
\end{align*}
or equivalently their continuous-time counterparts:
\begin{align*}
  A_K^{*} \bar{Q}_1 + \bar{Q}_1A_K & = - K^{*}RK \\
   A_K^{*} \bar{Q}_2 + \bar{Q}_2A_K & = - (C+DK)^{*}Q(C+DK),
\end{align*}
where $A_K := A + BK_{\Lambda}$. Note, however, that here it is required that $C$ is admissible for $A_K$ in order for the Lyapunov equations to have solutions \cite[Thm. 5.1.1]{TucWeiBook}.

Finally, the input and output terminal penalties are given by
$\left\langle x_{k+N},\right.$ $\left.\bar{Q}_1x_{k+N} \right\rangle$ and
$\left\langle  x_{k+N}, \bar{Q}_2 x_{k+N} \right\rangle$,
respectively. Thus, the quadratic formulation of the MPC problem is
given as in the stable case, except that in $H$ and $P$ the operator
$\bar{Q}$  must be replaced with $\bar{Q}_1 + \bar{Q}_2$.

Note that the full state feedback $u = Kx$ optimally solves the
\textit{unconstrained} minimization problem \eqref{eq:mpc}. Thus,
in order to utilize it in the constrained setting, we need to first
assume that the system is stabilizable by a sequence of inputs satisfying the input constraints. Under this assumption, MPC is utilized to steer the system into a region where $u_{\min} \leq Kx \leq u_{\max}$, at which point we can switch from MPC to the state feedback control. Switching to feedback control is in fact necessary as it has been shown in \cite[Thm. 1.3]{RebTow97} that an arbitrary stabilizable system cannot be stabilized by piecewise polynomial control. The existence of a constrained stabilizing input sequence can be guaranteed by allowing sufficiently high-gain inputs to cancel out the unstable dynamics of the system.

\begin{thm}
\label{thm:2}
Assume that the system \eqref{eq:linsys} is
stabilizable by a sequence of inputs satisfying the input constraints. Then, the dual-mode control consisting of MPC and optimal state feedback optimally stabilizes the system while satisfying the input constraints.
\begin{proof}
As the stabilization cost is included in the MPC problem, the optimal
solutions of \eqref{eq:jkfin} asymptotically steer the state of the system towards zero. Once the state reaches the region where state feedback satisfies the input constraints, MPC can be switched to unconstrained gain based controller to finalize stabilization.
\end{proof}
\end{thm}

In practice, finding the optimal feedback $K$ is rather challenging as the Riccati equation \eqref{eq:cric} can rarely be solved in analytic closed-form. Moreover, in general there is no guarantee that $C$ would be admissible for the stabilized semigroup generated by $A+BK_\Lambda$. However, some other stabilizing feedback can be used as a terminal penalty and stabilizing feedback as well, such as output feedback $u_k = K_yy_k$. This is a valid choice as regularity of the system is preserved under output feedback (see \cite{Wei94a}), and rather straightforward computations using Sherman-Morrison-Woodbury formula show that $A_d + B_dK_y(I - D_dK_y)^{-1}C_d$ corresponds to the Cayley-Tustin discretization of $A + BK_y(I-DK_y)^{-1}C$, i.e., $A$ after output feedback. Apart from optimality, the result of Theorem \ref{thm:2} holds for any stabilizing feedback.

\section{Wave Equation} \label{sec:waveeq}

As an example of a stable system, consider the wave equation on a 1-D spatial domain $\zeta \in [0, 1]$
with viscous damping at one end and boundary control $u$ and boundary
observation $y$ at the other end given by:
\begin{subequations}
  \label{eq:wavebcs}%
  \begin{align}
    \frac{\partial^2}{\partial t^2}w(\zeta,t)  & =  \frac{1}{\rho(\zeta)}\frac{\partial}{\partial\zeta}\left(T(\zeta)\frac{\partial}{\partial\zeta}w(\zeta,t) \right) \\
	0 & = T(\zeta)\frac{\partial}{\partial \zeta}w(1,t) + \frac{\kappa}{\rho}\frac{\partial}{\partial t}w(1,t) \\
    u(t) & = T(\zeta)\frac{\partial}{\partial\zeta}w(0,t) \\
    y(t) & = \frac{\partial}{\partial t}w(0,t), \label{eq:waveout}
\end{align}
\end{subequations}
where $\kappa > 0$. For simplicity we assume that the mass density
$\rho$ and the Young's modulus $T$ are constants. We further assume
that $\kappa \neq \sqrt{\rho T}$, which will be needed in Section \ref{waveeq:lyap}.

In order to write \eqref{eq:wavebcs} in the state-space form \eqref{eq:linsys}, let us
first define a new state variable $x = [x_1, \, x_2]^T := [\rho \partial_t w, \, \partial_\zeta w]^T$
with state space $X = L^2(0,1;\mathbb{R}^2)$ and an auxiliary matrix operator $\mathcal{H}(\zeta):= \operatorname{diag}(\rho(\zeta)^{-1}, \; T(\zeta))$.
Now the operator $A$ can be defined as:
$$
 Ax(\zeta, t) :=
 \begin{bmatrix}
 0 & 1 \\ 1 & 0
 \end{bmatrix}
  \frac{\partial}{\partial\zeta}\left( \mathcal{H}(\zeta)x(\zeta,t) \right)
$$
with domain
$$
\mathcal{D}(A) := \left\{x \in X: \mathcal{H}x
  \in H^1(0,1;\mathbb{R}^2),\;  x \in \mathcal{N}(\mathcal{B}) \right\},
$$
where $\mathcal{B}x := [Tx_2(1,\cdot) + \kappa\rho^{-1} x_1(1,\cdot), Tx_2(0,\cdot)]^T$. It follows from
\cite[Thm. III.2]{VilZwa09} that $A$ is the generator of an
exponentially stable $C_0$-semigroup. 

By \cite[Rem. 10.1.5]{TucWeiBook}, the control operator $B \in \mathcal{L}(U,X_{-1})$ can be found by solving the abstract elliptic equation 
\begin{equation}
\label{eq:TucWeiB}
\begin{split}
\rho^{-1}f'_1(\zeta) & = sf_2(\zeta), \quad Tf'_2(\zeta) = sf_1(\zeta)  \\
Tf_2(1) + \kappa\rho^{-1}f_1(1) & = 0, \qquad Tf_2(0) = 1,
\end{split}
\end{equation}
the solution of which satisfies $f = (s - A_{-1})^{-1}B$ for any $s \in \rho(A)$. While the abstract equation can be utilized in finding $B$, we note that by solving the equation for $s = \delta$ we directly obtain the discretized control operator as $B_d = \sqrt{2\delta}f$.

The observation operator $C \in \mathcal{L}(X_1, Y)$ is simply defined as $Cx := \rho^{-1}x_1(0)$ on $\mathcal{D}(A)$ and can be extended as such to $\mathcal{C}x := \rho^{-1}x_1(0)$ on $H^1(0,1;\mathbb{R}^2)$ (see \cite[Rem 4.11]{Wei94b}). Furthermore, the transfer function of \eqref{eq:wavebcs} is given by $G(s) = \mathcal{C}(s - A_{-1})^{-1}B$ which can be evaluated by solving \eqref{eq:TucWeiB} and applying $\mathcal{C}$ into the solution, that is $G(s) = \mathcal{C}f$. Equivalently by \cite[Def. 5.6, Thm. 5.8]{Wei94b}, the transfer function can be written using the $\Lambda$-extension of $C$ as $G(s) = C_\Lambda(s - A_{-1})^{-1}B + D$, where $D = \displaystyle \lim_{s\to \infty} G(s)$. We will utilize the former way of evaluating the transfer function for $D_d = G(\delta)$. Finally, note that no extension of $C$ is required for computing $C_d$ as $\mathcal{R}(\delta - A)^{-1} \subset \mathcal{D}(C)$.

\subsection{Discretized Operators} \label{waveeq:disc}

Assume that $\rho$ and $T$ are constants and consider the equation
$\dot{x}(t) = Ax(t)$. Using the  Laplace transform yields
$$
sx(\zeta,s) - x(\zeta,0) = \frac{\partial}{\partial \zeta}\left( 
\begin{bmatrix}
0 & T \\ \rho^{-1} & 0
\end{bmatrix}x(\zeta,s)
 \right),
$$
that is,
$$
\frac{\partial}{\partial \zeta}x(\zeta,s) = 
\begin{bmatrix}
0 & \rho s \\ T^{-1}s & 0
\end{bmatrix}x(\zeta,s) - 
\begin{bmatrix}
0 & \rho \\ T^{-1} & 0
\end{bmatrix}x(\zeta,0).
$$
The above is an ordinary differential equation of the form:
$$
  \frac{\partial}{\partial \zeta}x(\zeta,s) = \overline{A}x(\zeta,s) - \overline{B}x(\zeta,0),
$$
the solution of which is given by:
\begin{equation}
\label{eq:res1}
x(\zeta,s) = e^{\overline{A}\zeta}x(0,s) -
\int\limits_0^{\zeta}e^{\overline{A}(\zeta-\eta)} \overline{B}x(\eta,0)d\eta
\end{equation}
where:
$$
e^{\overline{A}\zeta} = 
\begin{bmatrix}
\cosh \left( \sqrt{\frac{\rho}{T}}s\zeta \right) & \sqrt{\rho T}\sinh
\left( \sqrt{\frac{\rho}{T}}s\zeta \right) \\
\left( \sqrt{\rho T} \right)^{-1}\sinh \left(
  \sqrt{\frac{\rho}{T}}s\zeta \right) & \cosh \left( \sqrt{\frac{\rho}{T}}s\zeta \right)
\end{bmatrix}.
$$
Recall that $\mathcal{D}(A)$ has the boundary conditions $Tx_2(1) +
\frac{\kappa}{\rho}x_1(1) = 0$ and $Tx_2(0) = 0$, based on which
$x(0,s)$ in \eqref{eq:res1} can be solved. Eventually, 
\eqref{eq:res1} is given by:
\begin{align*}
x(\zeta,s) & = \frac{\rho}{\sqrt{\rho T}\sinh \left(
               \sqrt{\frac{\rho}{T}}s \right) + \kappa\cosh\left(
               \sqrt{\frac{\rho}{T}}s \right)}
\begin{bmatrix}
  \cosh \left( \sqrt{\frac{\rho}{T}}s\zeta \right) \\
  \left( \sqrt{\rho T} \right)^{-1}\sinh \left( \sqrt{\frac{\rho}{T}}s\zeta \right)
\end{bmatrix}\times \\
&  \quad \int\limits_0^1 \textstyle \left( \frac{\kappa}{\sqrt{\rho T}}\sinh
                         \left(\sqrt{\frac{\rho}{T}}s(1-\eta)\right) +
                         \cosh\left(\sqrt{\frac{\rho}{T}}s(1-\eta)\right)\right)
                         x_1(\eta,0)\\
  & \textstyle \qquad + \left(
  \kappa\cosh\left(\sqrt{\frac{\rho}{T}}s(1-\eta)\right) + \sqrt{\rho
  T}\sinh\left(\sqrt{\frac{\rho}{T}}s(1-\eta)\right)
  \right)x_2(\eta, 0)d\eta \\
  & \quad - \int\limits_0^\zeta 
\begin{bmatrix}
	\sqrt{\frac{\rho}{T}}\sinh\left(\sqrt{\frac{\rho}{T}}s(\zeta-\eta)\right) & \rho\cosh\left( \sqrt{\frac{\rho}{T}}s(\zeta-\eta)\right) \\
	T^{-1}\cosh\left( \sqrt{\frac{\rho}{T}}s(\zeta-\eta)\right) & \sqrt{\frac{\rho}{T}}\sinh\left( \sqrt{\frac{\rho}{T}}s(\zeta-\eta)\right)
      \end{bmatrix} x(\eta,0)d\eta \\
      & := (s - A)^{-1}x(\zeta,0),
\end{align*}
which yields the expression for the resolvent operator, from  which we
also obtain the operator $A_d = -I + 2\delta(\delta - A)^{-1}$.

In order to find the expression for $B_d$, we recall that the solution $f$ of \eqref{eq:TucWeiB} satisfies $f = (s - A_{-1})^{-1}B$ for all $s \in \rho(A)$. Thus, setting $s = \delta$, we obtain that $B_d$ is given by
    \begin{align*}
      B_d & = \frac{-\sqrt{2\delta}}{\sqrt{\rho T}\sinh\left(
        \sqrt{\frac{\rho}{T}}\delta \right)+\kappa\cosh \left( \sqrt{\frac{\rho}{T}}\delta \right)} \times
      \\ & \quad
\begin{bmatrix}
\rho \cosh \left( \sqrt{\frac{\rho}{T}}\delta(\zeta-1) \right) + \kappa
\sqrt{\frac{\rho}{T}}\sinh \left( \sqrt{\frac{\rho}{T}}\delta(\zeta+1)
\right) \\
\sqrt{\frac{\rho}{T}}\sinh \left( \sqrt{\frac{\rho}{T}}\delta(\zeta-1)
\right) - \frac{\kappa}{T}\cosh \left( \sqrt{\frac{\rho}{T}}\delta(\zeta-1) \right)
\end{bmatrix}.
\end{align*}
Furthermore, we obtain:
\begin{align*}
C_dx(\zeta) & = \frac{\sqrt{2\delta}}{\sqrt{\rho T}\sinh \left(
    \sqrt{\frac{\rho}{T}\delta} \right) \kappa\cosh \left(
    \sqrt{\frac{\rho}{T}}\delta \right)}\times \\
  &  \quad \int\limits_0^1 \textstyle \left( \frac{\kappa}{\sqrt{\rho T}}\sinh
                         \left(\sqrt{\frac{\rho}{T}}\delta(1-\eta)\right) +
                         \cosh\left(\sqrt{\frac{\rho}{T}}\delta(1-\eta)\right)\right)
                         x_1(\eta)\\
  & \textstyle \qquad + \left(
  \kappa\cosh\left(\sqrt{\frac{\rho}{T}}\delta(1-\eta)\right) + \sqrt{\rho
  T}\sinh\left(\sqrt{\frac{\rho}{T}}\delta(1-\eta)\right)
  \right)x_2(\eta)d\eta.
\end{align*}
Finally, in order to compute $D_d = G(\delta)$ we evaluate the transfer function at $\delta$ by $\mathcal{C}(\delta - A_{-1})^{-1}B$, where $(\delta - A_{-1})^{-1}B$ is again obtained by solving \eqref{eq:TucWeiB} for $s = \delta$ (which is equal to $B_d/\sqrt{2\delta}$). Thus, $D_d$ is given by
\begin{equation}
  \label{eq:Ddisc}
  D_d := -\frac{1}{\sqrt{\rho T}}\frac{\kappa \sinh \left(
      \sqrt{\frac{\rho}{T}}\delta \right) + \sqrt{\rho T}\cosh \left(
      \sqrt{\frac{\rho}{T}}\delta \right)}{\sqrt{\rho T}\sinh \left(
      \sqrt{\frac{\rho}{T}}\delta \right) + \kappa\cosh \left( \sqrt{\frac{\rho}{T}}\delta \right)}.
\end{equation}

In order to verify that the wave system \eqref{eq:wavebcs} is regular, we first utilize the fact that the $C_0$-semigroup generated by $A$ is invertible \cite[Def. 2.7.1]{TucWeiBook}, which we will show in the next section when finding the spectrum of $A$. Thus, we can utilize \cite[Thm. 5.2.2, Cor. 5.2.4]{TucWeiBook} by which $B$ and $C$ are admissible if $\|(s-A_{-1})^{-1}B\|$ and $\|C(s-A)^{-1}\|$, respectively, are uniformly bounded along a vertical line on the complex right half plane, which are easy to verify based on the expressions derived for $B_d$ and $C_d$. Finally, we note 
that $\displaystyle \lim_{\delta\to\infty}G(\delta) = -(\rho T)^{-1/2}$
to conclude that \eqref{eq:wavebcs} indeed is a regular linear system.

The adjoints of the discretized operators can be computed based on the property that the adjoint $P^*$ of an operator $P$ satisfies $<Px,y> = <x,P^*y>$ with respect to the corresponding inner products. Here the state-space $X$ is equipped with the $L_2$ inner product, and the input and output spaces are equipped with the real scalar product. The computations for finding the adjoints are straightforward and will be omitted here for brevity. An explicit example of computing the adjoint operators can be found e.g. in \cite{HumDubACC18}.

\subsection{Solution of the Lyapunov equation} \label{waveeq:lyap}

In this section, we derive the positive solution for the continuous Lyapunov
equation \eqref{eq:lyacont},
which is realized by utilizing the spectral representation of $A$. Let us
at first find the eigenvalues and eigenvectors of the operator $A$. A
direct computation shows that the solution of the 
eigenvalue equation $A\phi_k = \lambda_k\phi_k$ is of the form:
\begin{align*}
\phi_{1,k}(\zeta) & = \alpha \exp \left( \sqrt{\frac{\rho}{T}}\lambda_k\zeta
  \right) + \beta\exp \left( \sqrt{\frac{\rho}{T}}\lambda_k\zeta
             \right) \\
  \phi_{2,k}(\zeta) & = \frac{\alpha}{\sqrt{\rho T}}\exp \left(
               \sqrt{\frac{\rho}{T}}\lambda_k\zeta \right) -
               \frac{\beta}{\sqrt{\rho T}}\exp \left( \sqrt{\frac{\rho}{T}}\lambda_k\zeta \right).
\end{align*}
Since $\phi_k\in \mathcal{D}(\mathcal{A})$, we must have
$\phi_{2,k}(0) = 0$, which yields $\alpha = \beta$. Thus, the
eigenvectors of $A$ are of the form:
$$
\phi_k(\zeta) = 
\begin{bmatrix}
  \cosh \left( \sqrt{\frac{\rho}{T}}\lambda_k\zeta \right) \\
  \frac{1}{\sqrt{\rho T}}\sinh \left( \sqrt{\frac{\rho}{T}}\lambda_k\zeta \right)
\end{bmatrix},
$$
and the eigenvalues $\lambda_k$ are determined from the condition
$T\phi_{2,k}(1) = -\frac{\kappa}{\rho}\phi_{1,k}(1)$, i.e.,
$$
\sqrt{\frac{T}{\rho}}\sinh \left( \sqrt{\frac{\rho}{T}}\lambda_k
\right) + \frac{\kappa}{\rho}\cosh \left(
  \sqrt{\frac{\rho}{T}}\lambda_k \right) = 0.
$$
Using the exponential form of the hyperbolic functions we obtain that
one of the eigenvalues is given by:
\begin{equation}
\label{eq:lam0}
\lambda_0  =\frac{1}{2}\sqrt{\frac{T}{\rho}}\log \left(
  \frac{\sqrt{\rho T} - \kappa}{\sqrt{\rho T}+ \kappa} \right),
\end{equation}
which is real if $\kappa < \sqrt{\rho T}$. Finally, by the periodicity of
the exponential function along the imaginary axis, we obtain that in
general the eigenvalues are given by
$\lambda_k = \lambda_0 + \sqrt{T/\rho}k\pi i$ for $k \in \mathbb{Z}$, which also implies that the semigroup generated by $A$ is invertible (see, e.g., \cite[Prop. 2.7.8]{TucWeiBook}) as stated in the previous section.

We note that damped wave equations have been considered, e.g., in \cite{CoxZua95} and \cite[Sect. 4]{XuCWei11} -- both referring to the original work by Rideau \cite{RidPhD} -- where similar spectra were obtained. Furthermore, it can be seen from \eqref{eq:lam0} that the assumption $\kappa\neq\sqrt{\rho T}$ is required to ensure $\sigma(A)\neq \emptyset$, which is further required by \cite[Thm. 3.5]{CoxZua95} to ensure that the eigenvectors of $A$ constitute a Riesz basis for $X$. Indeed, we can define an invertible operator:
$$
M := 
\begin{bmatrix}
\cosh \left( \sqrt{\frac{\rho}{T}}\lambda_0\zeta \right) & -\sqrt{\rho
T}\sinh \left( \sqrt{\frac{\rho}{T}}\lambda_0\zeta \right) \\ i\sinh
\left( \sqrt{\frac{\rho}{T}}\lambda_0\zeta \right) & -i \sqrt{\rho
  T}\cosh\left( \sqrt{\frac{\rho}{T}}\lambda_0\zeta \right) 
\end{bmatrix},
$$
so that $M\phi_k = [ \cos (k\pi\zeta), \;  \sin (k\pi\zeta)]$
is an orthonormal basis in $X$, and the biorthogonal sequence
\cite[Def. 2.5.1]{TucWeiBook} $\left(
  \bar{\phi}_{k} \right)$ to $\left( \phi_k \right)$ is given by $\bar{\phi}_{k}
= M^{*}M\phi_{k}$.

Let us now return to the Lyapunov equation and apply it to an
arbitrary $x \in \mathcal{D}(A)$: 
$$
A^{*} \bar{Q}x + \bar{Q}Ax + C^{*}QCx = 0.
$$
By \cite[Prop. 2.5.2]{TucWeiBook}, we can write $x = \sum_{k\in \mathbb{Z}} \left\langle z, \bar{\phi}_k \right\rangle\phi_k$ for every $x \in X$, which yields:
$$
\sum_{k\in \mathbb{Z}} \left( A^{*}\bar{Q}  \left\langle x, \bar{\phi}_k
  \right\rangle\phi_k + \bar{Q}A \left\langle x, \bar{\phi}_k
  \right\rangle\phi_k + C^{*}Q C \left\langle x,
    \bar{\phi}_k \right\rangle\phi_k \right) = 0,
$$
which by utilizing \cite[Prop. 2.6.3]{TucWeiBook} further yields:
$$
\sum_{k\in \mathbb{Z}} \left( (A^{*} + \lambda_{k})\bar{Q}  \left\langle x, \bar{\phi}_k
  \right\rangle\phi_k + C^{*}Q C \left\langle x,
    \bar{\phi}_k \right\rangle\phi_k \right) = 0.
$$
The above especially holds if $(A^{*} + \lambda_k)\bar{Q}\left\langle x, \bar{\phi}_k
  \right\rangle\phi_k = -C^{*}QC\left\langle x, \bar{\phi}_k
  \right\rangle\phi_k$ for all $k \in \mathbb{Z}$. Thus, for an
  arbitrary $k \in \mathbb{Z}$, we obtain:
$$
\bar{Q}\left\langle x, \bar{\phi}_k
  \right\rangle\phi_k = (-\lambda_k - A^{*})^{-1}C^{*}QC\left\langle x, \bar{\phi}_k
  \right\rangle\phi_k.
$$
As $A$ is densely defined and $-\bar{\lambda}_{k} \in \rho(A)$ since
$\lambda_k \in \sigma(A)$, we have by \cite[Prop. 2.8.4]{TucWeiBook}
that $(-\lambda_k - A^{*})^{-1} = \left((- \bar{\lambda}_k -
  A)^{-1}\right)^{*}$, so we obtain :
\begin{align*}
\bar{Q}\left\langle x, \bar{\phi}_k
  \right\rangle\phi_k & = \left( (- \bar{\lambda}_k - A)^{-1}
  \right)^{*} C^{*}QC\left\langle x, \bar{\phi}_k
  \right\rangle\phi_k \\ & = \left( C(-\bar{\lambda}_k - A)^{-1} \right)^{*}QC\left\langle x, \bar{\phi}_k
  \right\rangle\phi_k.
\end{align*}
Finally, summation over $k\in \mathbb{Z}$ yields the solution:
\begin{equation}
\label{eq:lyapsol}
\bar{Q}x = \sum_{k\in \mathbb{Z}} \left\langle x, \bar{\phi}_k \right\rangle\left( C(-\bar{\lambda}_k - A)^{-1} \right)^{*}QC\phi_k.
\end{equation}
Note that as $C\phi_k = 1$ and $C(-\lambda_k^* -
A)^{-1}$ is uniformly bounded for all $k\in\mathbb{Z}$, the series in
\eqref{eq:lyapsol} is convergent. Thus, for any $x\in X$ we may approximate:
$$
\bar{Q}x \approx \bar{Q}_Mx := \sum_{k=-M}^M \left\langle x, \bar{\phi}_k \right\rangle\left( C(-\bar{\lambda}_k - A)^{-1} \right)^{*}QC\phi_k,
$$
and it holds that $\displaystyle
\lim_{M\to\infty}\|\bar{Q}x - \bar{Q}_Mx\|_{L_2}=0$, by which we can
evaluate \eqref{eq:lyapsol} to an arbitrary precision $\epsilon > 0$ by
choosing a sufficiently large $M$. A suitable value for $M$ can
determined, e.g., by numerical experiments. 

\subsection{Simulation results for the wave equation} \label{waveeq:sim}

Consider the wave equation \eqref{eq:wavebcs}	 with the parameter choices $\rho = T = 1$ and $\kappa = 0.75$. For the MPC, choose the optimization horizon as $N=15$ and choose the input and output weights as $R=10$ and $Q=0.5$, respectively. For the Cayley-Tustin discretization, choose $h=0.075$ so that $\delta \approx 26.67$. For numerical integration, an adaptive approximation of $d\zeta$ is used with $519$ nodal points. To approximate the solution of the Lyapunov equation \eqref{eq:lyapsol}, we choose $M = 100$. The initial conditions for the wave equation are given by $\partial_t w(\zeta) = \cos(\pi\zeta)$ and  $\partial_\zeta w(\zeta) = \sin(\frac{1}{2}\pi\zeta)$.

\begin{figure}[htbp]
\begin{center}
\includegraphics[width=\columnwidth]{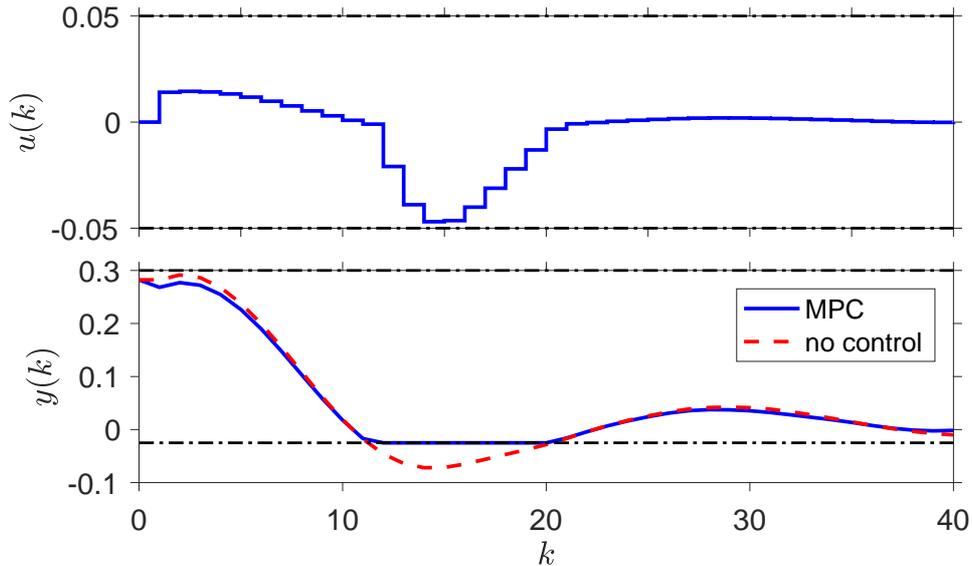}
\caption{Above: MPC inputs $u(k)$ and the input
constraints. Below: MPC and uncontrolled outputs and the output constraints.}
\label{fig:waveuy}
\end{center}
\end{figure}

The input and output constraints $-0.05 \leq u_k \leq 0.05$ and $-0.025 \leq y_k \leq 0.3$ are displayed in Figure \ref{fig:waveuy} along with the control inputs $u(k)$ obtained from the MPC problem. The outputs of the system under the MPC and under no control are displayed as well. It can be seen that the MPC makes the output decay slightly faster in the beginning. Then control is imposed to satisfy the output constraints while the uncontrolled output violates them. Finally, a minor stabilizing control effort is imposed before both the MPC input and the output decay to zero. Naturally the uncontrolled output decays to zero as well due to the exponential stability of the considered system.

\begin{figure}[htbp]
\begin{center}
\includegraphics[width=\columnwidth]{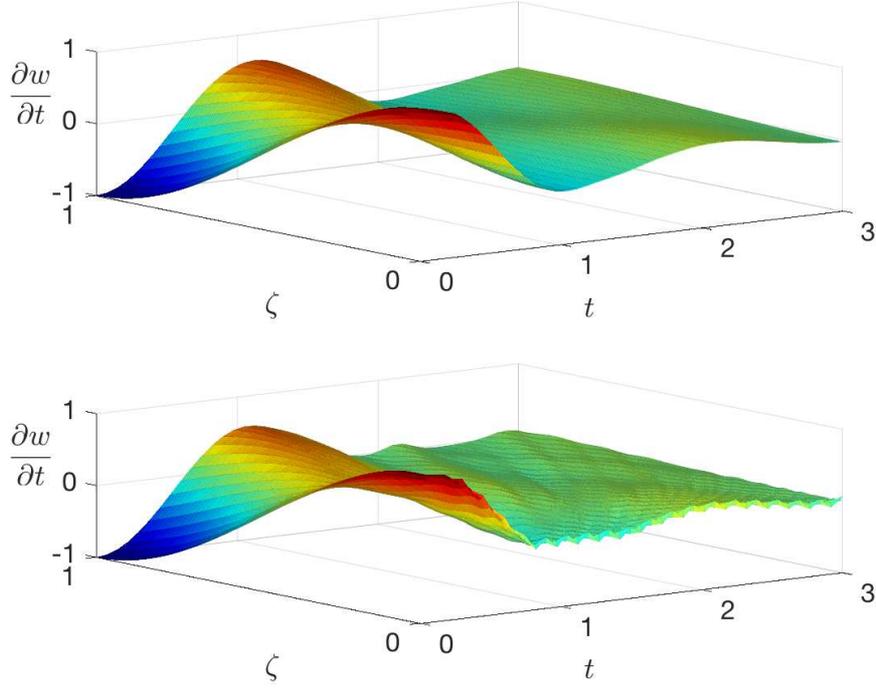}
\caption{Above: the velocity profile of the wave equation without control. Below: the velocity profile under the model predicting control law.}
\label{fig:waves}
\end{center}
\end{figure}

Figure \ref{fig:waves} displays the velocity profiles of the system under the model predictive control law and without control. No substantial differences can be observed in the velocity profiles, which is rather expected as the outputs in Figure \ref{fig:waveuy} were rather close to one another. We note that under the nonsmooth MPC inputs the wave system \eqref{eq:wavebcs} only admits a mild solution (a solution on $X_{-1}$) \cite[Prop. 4.2.5]{TucWeiBook}, due to which oscillations can be seen in the controlled velocity profile below.

\section{Tubular reactor with recycle} \label{sec:tub}

As another example, consider a model of the mixing of flow in a tubular reactor:
\begin{subequations}
\label{eq:tub}%
\begin{align}
\frac{\partial}{\partial t}x(\zeta,t) & = -v
                                        \frac{\partial}{\partial\zeta}x(\zeta,t)
                                        + \alpha x(\zeta,t) \\
  x(0,t) & = rx(1,t) + (1-r)u(t) \label{eq:tubu}\\ 
  y(t) & = x(1,t)
\end{align}
\end{subequations}
on $X = L^2(0,1;\mathbb{R})$, where $v, \alpha > 0$ are velocity and reaction terms, respectively, and $0 < r < 1$ is a recycle term. The system has its spectrum at $\{\alpha + v\log(r) + v2n\pi i\}_{n\in\mathbb{Z}}$, which for parameters $v = 1, \alpha = 1/2$ and $r = 2/3$ is on the right half plane. Thus, the system is unstable but exponentially stabilizable, e.g., by
output feedback $u(t) = -y(t)$. Under this feedback, \eqref{eq:tubu}
changes to $x(0,t) = (2r-1)x(1,t)$ but otherwise the system remains
the same. Note that the recycle term $r$ has the largest effect on the location of the spectrum as in the extremal case $r = 0$ the spectrum would be empty.

Similar to the wave equation in Section \ref{sec:waveeq}, we can
compute the resolvent operator and find the discretized operators
$(A_d,B_d,C_d,D_d)$ and their adjoints. Since output feedback is used
as a stabilizing terminal cost and in this case $D=0$, for the terminal penalty one needs to solve the Lyapunov
equation $A_s^{*} \bar{Q} + \bar{Q}A_s = - C^{*}(Q+R)C$,
where $A_s$ is the generator of the exponentially stable
$C_0$-semigroup corresponding to the boundary control system
\eqref{eq:tub} under output feedback $u(t) = -y(t)$. This can be done
as in Section \ref{waveeq:lyap}, except that the
normalized eigenvectors of $A_s$ already form an orthonormal basis
in $X = L^2(0,1;\mathbb{R})$.

Similar to the wave equation example \eqref{eq:wavebcs}, the eigenvalues of \eqref{eq:tub} are located (before and after stabilization) on a vertical line on the complex plane, and thus, the semigroup associated with the system is invertible \cite[Prop. 2.7.8]{TucWeiBook}. Thus, the admissibility of the operators $B$ and $C$ can be determined here as well based on $B_d$ and $C_d$, which together with the fact that $\displaystyle \lim_{s\to\infty}G(s) = 0$ concludes that the system \eqref{eq:tub} is regular.

For the MPC problem formulation, the weights are chosen as $Q = 2$ and
$R = 10$, and the input constraints are
given by $-0.15 \leq u_k \leq 0.05$ while no output constraints are imposed. The
optimization horizon is chosen as $N = 10$, and for approximation of
the solution of the Lyapunov equation, $201$ eigenvectors of $A_s$ are
used. For the Cayley-Tustin discretization, we choose $h = 0.1$ so that
$\delta = 20$. The initial condition is given by $x_0(\zeta) = \frac{1}{2}\sin(\pi\zeta)$. For numerical integration, an adaptive approximation of
$d\zeta$ is used with $510$ nodal points.

\begin{figure}[htbp]
\begin{center}
\includegraphics[width=\columnwidth]{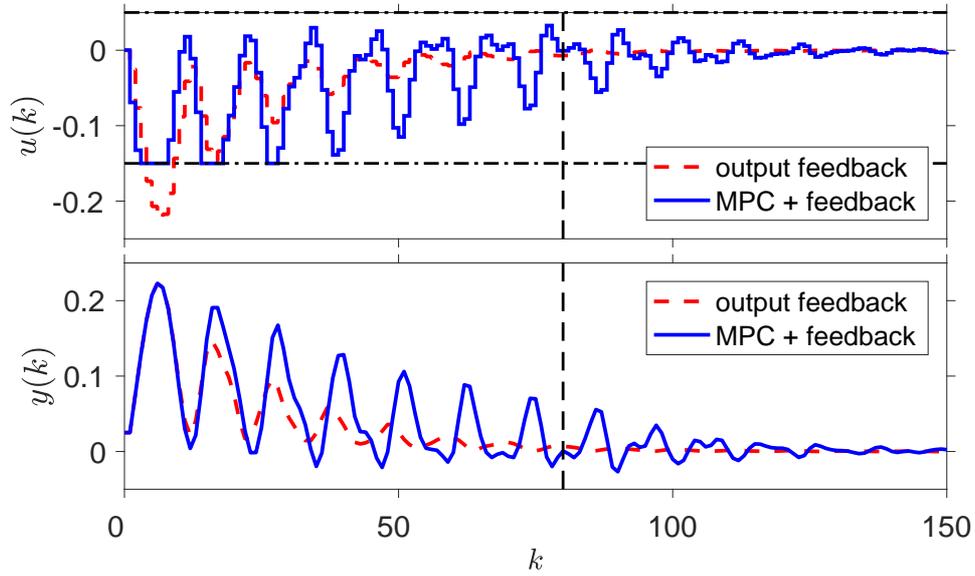}
\caption{Above: dual-mode inputs, the input
  constraints and the output feedback. Below: outputs
  of the system under the dual-mode control and output feedback.}
\label{fig:uytub}
\end{center}
\end{figure}

In Figure \ref{fig:uytub}, the dual-mode inputs and the outputs of the system
under the dual-mode control are presented. For comparison, the output
feedback control and the output under the feedback control are also
presented. It can be seen that while the output feedback stabilizes
the system faster, it does not satisfy the input constraints early on
in the simulation. In the dual-model control, the MPC inputs first steer the output close to
zero while satisfying the input constraints, and then at $k=80$ it is switched to output feedback $u = -y$ which completes the stabilization.

\begin{figure}[htbp]
\begin{center}
\includegraphics[width=\columnwidth]{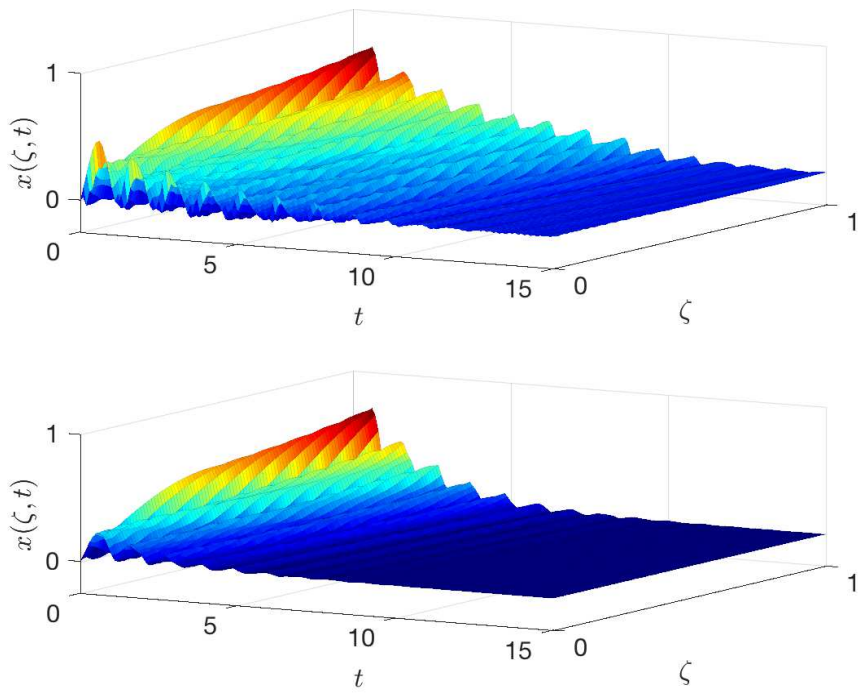}
\caption{Above: the state profile of the tubular reactor under the-dual mode control. Below: the state profile under the output feedback.}
\label{fig:tubs}
\end{center}
\end{figure}

In Figure \ref{fig:tubs}, the state profiles
of the tubular reactor are displayed under the dual-mode and the feedback controls. The states behave according to what could be expected based on the outputs, that is, both states decay asymptotically to zero and the state under output feedback decays faster. Similar to the wave system, the nonsmooth inputs admit the system to only have a mild solution, due to which there are oscillations in the state profiles below.

\section{Conclusions} \label{sec:concl}

In this work, a linear model predictive controller for regular linear
systems was designed, and it was shown that for stable systems, stability
of the zero output regulator follows from the finite-dimensional
MPC theory. For stabilizable systems, constrained stabilization was achieved by dual-mode control consisting of MPC and stabilizing feedback. The MPC design was demonstrated on an illustrative example where it was implemented for the boundary controlled wave equation. Constrained stabilization was demonstrated on a tubular reactor which had solely unstable eigenvalues. The performances of the control strategies were illustrated with numerical simulations.

\bibliographystyle{plain}

\end{document}